\documentclass[a4paper,11pt]{article}

\usepackage{amsmath}
\usepackage{amsthm}
\usepackage{amssymb}
\usepackage{graphics}
\usepackage{color}
\usepackage{epsfig}

\DeclareMathOperator{\link}{link}

\DeclareMathOperator{\pure}{pure}

\newcommand{\sm}{\!\setminus\!}

\theoremstyle{plain}
\newtheorem{thm}{Theorem}[section]
\newtheorem{cor}[thm]{Corollary}
\newtheorem{prop}[thm]{Proposition}
\newtheorem{lemma}[thm]{Lemma}
\theoremstyle{definition}
\newtheorem{defn}[thm]{Definition}

\newtheorem{question}[thm]{Question}
\theoremstyle{remark}
\newtheorem*{remark}{Remark}

\newcommand{\vanish}[1]{}

\title{ Obstructions to shellability, partitionability,  \\
        and sequential Cohen-Macaulayness}
\author{Masahiro Hachimori \\
        \small
        Graduate School of Systems and \\
        \small
        Information Engineering, \\
        \small
        University of Tsukuba \\
        \small
        Tsukuba, Ibaraki 305-8573, Japan \\
        \small
        {\tt hachi@sk.tsukuba.ac.jp}  
        \and
        Kenji Kashiwabara \\
        \small
        Department of General Systems Studies, \\
        \small
        University of Tokyo \\
        \small
        Komaba, Meguro, Tokyo 153-8902, Japan \\
        \small
        {\tt kashiwa@idea.c.u-tokyo.ac.jp}  
       }
\date{Jan 2011}

\begin{document}
\maketitle

\begin{abstract}
For a property $\cal P$ of simplicial complexes,
a simplicial complex $\Gamma$ is an obstruction to $\cal P$ if 
$\Gamma$ itself does not satisfy $\cal P$ but all of its proper restrictions
satisfy $\cal P$.
In this paper, we determine all obstructions to shellability of 
dimension $\le 2$, refining the previous work by Wachs.
As a consequence we obtain that
the set of obstructions to shellability,
that to partitionability and that to sequential Cohen-Macaulayness
all coincide for dimensions $\le 2$.
We also show that these three sets of obstructions coincide 
in the class of flag complexes.
These results show that the three properties, hereditary-shellability, 
hereditary-partitionability,
and hereditary-sequential Cohen-Macaulayness are equivalent for these classes.
\end{abstract}




\sloppy

\section{Introduction}
In this paper we treat a (finite) simplicial complex 
as an abstract simplicial complex, i.e.,
a set family $\Gamma$ over a finite set $E$ 
closed under taking subsets.
A member $\sigma$ of $\Gamma$ is a \textit{face} of $\Gamma$, whose
dimension $\dim \sigma$ is $\lvert\sigma\rvert - 1$.
A maximal face with respect to inclusion is a \textit{facet},
a $0$-dimensional face is a \textit{vertex}, 
and a $1$-dimensional face is an \textit{edge}.
For a vertex $\{v\}$, we also refer $v$ as a vertex for convenience 
if there will be no confusion,
and use the notation $V(\Gamma)$ for the set of all vertices of $\Gamma$.
(That is, $V(\Gamma) = \{v\in E: \{v\}\in \Gamma\}$.)
Note that the empty set $\emptyset$ 
is always a face of a simplicial complex, with $\dim\emptyset = -1$.
We sometimes denote a $k$-dimensional facet as a \textit{$k$-facet} for short.
The dimension $\dim\Gamma$ of $\Gamma$ is the maximum dimension of 
its faces.
$\Gamma$ is \textit{pure} if all its facets have dimension $\dim\Gamma$.
The \textit{pure $i$-skeleton} $\pure_i(\Gamma)$ of $\Gamma$
is the set of subfaces of all $i$-dimensional faces of $\Gamma$.
Note that the pure $i$-skeletons of $\Gamma$ are $i$-dimensional pure simplicial complexes.
For a face $\tau$ of $\Gamma$, the \textit{link} of $\tau$ in $\Gamma$ is
defined by $\link_\Gamma(\tau) = \{\eta\in\Gamma : \eta\cap\tau= \emptyset,
\eta\cup\tau\in\Gamma\}$.

In a simplicial complex, 
a sequence $\sigma_1, \sigma_2, \dotsc, \sigma_t$ of the facets
is called a \textit{shelling} if it satisfies that
$(\bigcup_{i=1}^{j-1}\overline{\sigma_i})\cap \overline{\sigma_j}$ is
a pure $(\dim \sigma_j - 1)$-dimensional subcomplex for every $2\le j\le t$,
where $\overline{\sigma}$ denotes the set of all faces of the face $\sigma$,
and a simplicial complex is \textit{shellable} if it has a shelling.
Note that this definition,
introduced by Bj\"orner and Wachs in \cite{B-W:1996, B-W:1996b},
applies both pure and nonpure cases,
extending the classical definition for pure complexes.
For the details of (nonpure) shellability together with 
additional terminology on simplicial complexes,
see \cite{B-W:1996}.

In the following, 
for $W\subseteq V(\Gamma)$, 
$\Gamma[W]$ denotes the \textit{restriction} of $\Gamma$ to $W$, 
i.e., the subcomplex on the vertex set $W$ that consists of the faces of $\Gamma$
contained in $W$.
For $U = V(\Gamma)\sm W$, we also use the notation 
$\Gamma\sm U =\Gamma[W]$,
the \textit{deletion} of $U$ from $\Gamma$.
Note that restricting and taking link are commutative, that is, 
we have $\link_{\Gamma}(\tau)[W] = \link_{\Gamma[W]}(\tau)$ for $\tau\in\Gamma[W]$.

Wachs introduced the concept of obstructions to shellability
in \cite{W:2000}.
A simplicial complex $\Gamma$ is an \textit{obstruction to shellability} 
if $\Gamma$ is not shellable but $\Gamma[W]$ is shellable
for all $W\subsetneq V(\Gamma)$.
(Equivalently, $\Gamma$ is an obstruction to shellability
if $\Gamma$ is not shellable and $\Gamma\sm U$ is shellable
for all $U\subseteq V(\Gamma)$ with $U\neq \emptyset$.)
In other words, an obstruction to shellability is a minimal simplicial complex
with respect to restriction that is nonshellable.
(The terminology ``\textit{minimally nonshellable complexes}'' is used
for obstructions to shellability in \cite{W:2007}.)
In \cite{W:2000}, Wachs stated the following.

\begin{thm}(\cite{W:2000})\mbox{} \label{thm:wachs}
\begin{enumerate}
\setlength{\itemsep}{0pt}
\item There are no obstructions to shellability of dimension 0.
\item The unique obstruction to shellability of dimension 1 is the simplicial
      complex on $4$ vertices with two disjoint 1-facets.
      (We denote this $1$-dimensional obstruction as $2K_2$.)
\item Obstructions to shellability exist for each dimension $\ge 1$.
\item Obstructions to shellability of dimension $2$ have at most $7$ vertices.
      Thus there are only finitely many obstructions to shellability
      of dimension $2$.
\end{enumerate}
\end{thm}

This result (iv) (with (ii)) of Wachs is important because it shows 
that the number of
obstructions to shellability is finite in the class of dimension $\le 2$.
But, on the other hand, 
her proof does not explicitly enumerate the individual obstructions
to shellability of dimension $2$.
Though the number of vertices is bounded by $7$,
it is not easy to determine all the obstructions by brute force checking.
In Section~\ref{sec:ots} of this paper, we determine all obstructions to shellability 
of dimension $\le 2$.
This result provides us a way to discuss further properties of obstructions
to shellability.

The concept of obstructions naturally applies to other properties of
simplicial complexes.
A simplicial complex $\Gamma$ is an \textit{obstruction to a property $\cal P$}
if $\Gamma$ does not satisfy $\cal P$ but $\Gamma[W]$ satisfies $\cal P$
for all $W\subsetneq V(\Gamma)$.
The obstructions to purity are already discussed by Wachs~\cite{W:2000}.
In Section~\ref{sec:otp-otCM}, 
we show that both
the set of obstructions to partitionability
and the set of obstructions to sequential Cohen-Macaulayness are 
the same as the set of obstructions to shellability
for dimensions $\le 2$,
as an application of our classification of 
obstructions to shellability of dimension $\le 2$.

A simplicial complex is called a \textit{flag complex} 
if all the minimal nonfaces are of size 2.
The class of flag complexes, which sometimes appear as clique complexes or
independence complexes of graphs,
has been attracting interests of many researchers coming from commutative algebra
as well as from combinatorics.
Recently, Woodroofe~\cite{W:2009} determined 
the obstructions to shellability 
in the class of flag complexes.
As a corollary of his result,
our discussion in Section~\ref{sec:otp-otCM} shows that
the set of obstructions to shellability, partitionability, and
sequentially Cohen-Macaulayness also coincide in the class of flag complexes.

For a property $\cal P$ of simplicial complexes, we say
a simplicial complex $\Gamma$ is \textit{hereditary-$\cal P$} if
all restrictions of $\Gamma$ (including $\Gamma$ itself) satisfy $\cal P$.
As is discussed in Section~\ref{sec:otp-otCM},
these results on obstructions show that
three properties, hereditary-shellability, hereditary-partitionability, and
hereditary-sequential Cohen-Macaulayness, are all equivalent
in the class of simplicial complexes of dimension at most $2$,
or in the class of flag complexes.

Further discussions toward higher dimensions are added in 
Section~\ref{sec:s-obstruction}.

\section{2-dimensional obstructions to shellability}
\label{sec:ots}

First we recall the following well-known but important
propositions about shellability.

\begin{prop}(\cite[Proposition~10.14]{B-W:1996})\label{prop:link}
Let $\Gamma$ be a shellable simplicial complex.
Then $\link_\Gamma(\tau)$ is shellable for any face $\tau\in\Gamma$.
\end{prop}

\begin{prop}(Rearrangement lemma~\cite[Lemma 2.6]{B-W:1996})\label{prop:rearrange}
Let $\Gamma$ be a shellable simplicial complex.
Then $\Gamma$ has a shelling $\sigma_1,\sigma_2,\dotsc,\sigma_t$
such that $\dim\sigma_i\ge\dim\sigma_j$ for $i\le j$.
\end{prop}

\begin{prop}(\cite[Theorem~2.9]{B-W:1996})\label{prop:purepart}
Let $\Gamma$ be a shellable simplicial complex.
Then $\pure_i(\Gamma)$ is shellable for all $0\le i\le \dim\Gamma$.
\end{prop}

\begin{prop}(\cite[Section~11.5]{B:1995}) \label{prop:homology}
Let $\Gamma$ be a pure shellable simplicial complex.
Then, for any $\tau\in\Gamma$,
$\widetilde{H}_k(\link_\Gamma(\tau))=0$ unless $k=\dim\link_\Gamma(\tau)$,
where $\widetilde{H}_k$ is the $k$-th reduced homology group 
(over $\mathbb{Z}$).
Especially, $\widetilde{H}_k(\Gamma)=0$ unless $k=\dim\Gamma$.
(In other words, a pure shellable simplicial complex is Cohen-Macaulay,
see Definition~\ref{defn:sCM}.)
\end{prop}

Also the following are easy to observe from the definition of shellability.

\begin{prop}
Every $0$-dimensional complex is shellable.
A pure $1$-dimensional complex is shellable if and only if it is connected.
\end{prop}

\begin{prop}\label{prop:1dim}
A $1$-dimensional complex is shellable if and only if 
$\pure_1(\Gamma)$ is connected.
\end{prop}

\begin{prop}\label{prop:2dim}
A $2$-dimensional complex $\Gamma$ is shellable if and only if 
$\pure_2(\Gamma)$ is shellable and $\pure_1(\Gamma)$ is connected.
\end{prop}

\begin{figure}
\begin{center}
\scalebox{.5}{\includegraphics{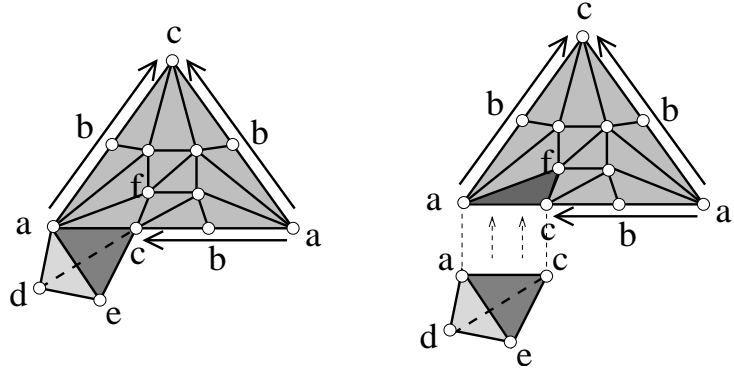}}
\caption{A nonshellable complex whose pure skeletons are all shellable.}
\label{fig:example}
\end{center}
\end{figure}

Remark that 
Propositions~\ref{prop:1dim} and \ref{prop:2dim} imply that the converse
of Proposition~\ref{prop:purepart} holds for dimensions $\le 2$, but
that this is not the case for dimensions $\ge 3$.
For example, Figure~\ref{fig:example} shows an example of 
a $3$-dimensional nonshellable complex whose pure skeletons are all shellable.
In the figure, $\{a,c,d,e\}$ is a $3$-dimensional facet, and the vertices
and edges are identified as indicated.
The upper part consisting of thirteen $2$-dimensional facets is 
shellable, but its shelling always has the facet $\{a,c,f\}$
(the dark facet in the right figure) in the last.
(This gadget is used in~\cite[Example 7]{H:2008} and the explanation 
of this property is essentially given there.)
But, the existence of the $3$-dimensional facet
enforces the shelling of the whole complex to start with 
the facet $\{a,c,d,e\}$
and the second facet should be $\{a,c,f\}$.
Thus the shelling cannot be completed.
On the other hand, it is not difficult to check that the pure skeletons
of the complex are all shellable.

From here we start the classification of 
2-dimensional obstructions to shellability.
First, we observe the following two lemmas.

\begin{lemma}\label{lemma:connected}
If $\Gamma$ is an obstruction to shellability, 
then it contains no isolated vertices.
Further, 
if $\dim\Gamma\ge 2$, then it is connected.
\end{lemma}
\begin{proof}
If $\Gamma$ has an isolated vertex $v$, then
it is easy to see that $\Gamma\sm\{v\}$ is shellable if and only if
$\Gamma$ is shellable.
This implies that 
an obstruction to shellability cannot have isolated vertices.

If $\dim\Gamma\ge 2$ and it is disconnected, 
we can choose two edges $\{u,v\}$ and $\{w,z\}$ belonging to different 
components since $\Gamma$ has no isolated vertices. 
This implies that  $\Gamma[\{u,v,w,z\}]$ is nonshellable and $1$-dimensional,
contradicting that $\Gamma$ is an obstruction to shellability of
dimension at least $2$.
\end{proof}

\begin{lemma}\label{lemma:2facets}
If $\Gamma$ is a $2$-dimensional obstruction to shellability,
then it has at least two $2$-facets.
\end{lemma}
\begin{proof}
By Lemma~\ref{lemma:connected}, $\Gamma$ is connected and thus
$\pure_1(\Gamma)$ is also connected.
By Proposition~\ref{prop:2dim}, this implies that $\pure_2(\Gamma)$
is nonshellable. For this, $\Gamma$ needs at least two 2-facets.
\end{proof}

The following proposition is important in our discussion of this section.

\begin{prop} \label{prop:edgeminimal}
Let $\Gamma$ be a $2$-dimensional obstruction to shellability,
and assume that $\Gamma$ has two vertices $v$ and $w$
such that $\{v, w\}$ is not an edge of $\Gamma$.
Then $\Gamma'=\Gamma \cup \{\{v,w\}\}$ is also an obstruction to shellability.
\end{prop}
\begin{proof}
$\Gamma$ is connected by Lemma~\ref{lemma:connected}, 
thus $\pure_1(\Gamma)$ is also connected.
Since $\Gamma$ is nonshellable, we have $\pure_2(\Gamma)$ nonshellable
by Proposition~\ref{prop:2dim}.
This implies $\Gamma'$ is nonshellable because
$\pure_2(\Gamma')=\pure_2(\Gamma)$.
On the other hand, we have $\Gamma[W]$ shellable for $W\subsetneq V(\Gamma)$,
hence $\pure_2(\Gamma[W])$ is shellable 
and $\pure_1(\Gamma[W])$ is connected.
This implies that $\pure_2(\Gamma'[W])$ is shellable since 
$\pure_2(\Gamma'[W])=\pure_2(\Gamma[W])$, 
and $\pure_1(\Gamma'[W])$ is connected because adding an edge 
between vertices does not make a connected graph into disconnected.
Thus we conclude that $\Gamma'[W]$ is shellable for $W\subsetneq V(\Gamma)$.
\end{proof}

We define a 2-dimensional obstruction to shellability to be \textit{edge-minimal} 
if it has a minimal set of edges with respect to inclusion relation.
By Proposition~\ref{prop:edgeminimal} above, all the obstructions to 
shellability of dimension 2 are determined if the edge-minimal obstructions
are specified.
We state our main theorem of this section
that gives the complete list of the
2-dimensional edge-minimal obstructions to shellability as follows.

\begin{figure}
\begin{center}
\scalebox{.4}{\includegraphics{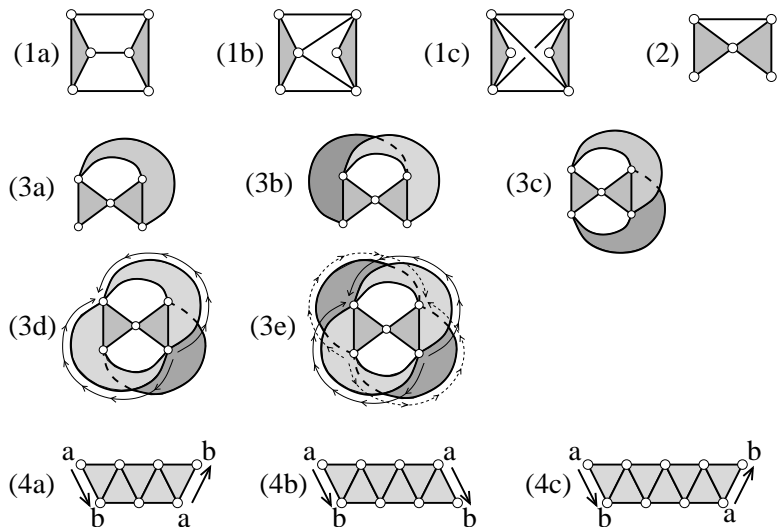}}
\end{center}
\caption{2-dimensional edge-minimal obstructions to shellability.
(In the figure, 
the pairs of edges indicated by the arrows are identified.)}
\label{fig:obstructions}
\end{figure}

\begin{thm} \label{thm:ots}
The complete list of $2$-dimensional edge-minimal obstructions to shellability
is given in Figure~\ref{fig:obstructions}.
Hence the obstructions to shellability of dimension $\le 2$
are the $1$-dimensional obstruction $2K_2$ of Theorem~\ref{thm:wachs} (ii)
and the simplicial complexes obtained by adding 
zero or more
edges to one of 
the complexes in Figure~\ref{fig:obstructions}.
\end{thm}

In the rest of this section we prove this theorem using several lemmas.

Let $\Gamma$ be an edge-minimal obstruction to shellability of dimension 2.
In a 2-dimensional simplicial complex,
we call an edge a \textit{nonboundary edge}
if it belongs to two or more 2-facets,
and a \textit{boundary edge} otherwise.
(We remark that we will also refer to an edge that does not belong to any 2-facets
as a boundary edge.)
We first discuss the case where
all the edges of $\Gamma$ are boundary edges in Lemma~\ref{lemma:boundarycase}.

\begin{lemma} \label{lemma:boundarycase}
Let $\Gamma$ be a $2$-dimensional edge-minimal obstruction to shellability 
with no nonboundary edges. 
Then $\Gamma$ is one of (1a), (1b), (1c), or (2) of Figure~\ref{fig:obstructions}.
\end{lemma}
\begin{proof}
By Lemma~\ref{lemma:2facets}, $\Gamma$ has at least two 2-facets $\sigma_1$ and $\sigma_2$.
From the fact that $\Gamma$ contains no nonboundary edges
it is easy to observe that $\pure_2(\Gamma[\sigma_1\cup\sigma_2])$ is nonshellable,
thus $\Gamma[\sigma_1\cup\sigma_2]$ is nonshellable by Proposition~\ref{prop:2dim}.
By the assumption that $\Gamma$ is an obstruction to shellability, 
there cannot exist a vertex outside of $\sigma_1 \cup \sigma_2$.

Assume that $\Gamma$ contains a third 2-facet $\sigma_3$.
Because $\Gamma$ has no nonboundary edges, any two 2-facets of $\Gamma$ cannot share
more than one vertices.
This implies that there exists a vertex of $\sigma_3$ that is not contained
in $\sigma_1\cup\sigma_2$, a contradiction.
Thus $\Gamma$ has exactly two 2-facets, and
there are no other vertices than those belonging to the two 2-facets.
We have the following two cases:
[Case 1] $V(\Gamma) = \{a,b,c,d,e,f\}$ 
and the 2-facets are $\sigma_1=\{a,b,c\}$ and $\sigma_2=\{d,e,f\}$,
or 
[Case 2]
$V(\Gamma) = \{a,b,c,d,e\}$ 
and the 2-facets are $\tau_1=\{a,b,c\}$ and $\tau_2=\{a,d,e\}$.

In both cases, $\pure_2(\Gamma)$ is nonshellable, thus $\Gamma$ is nonshellable
for any addition of edges by Proposition~\ref{prop:2dim}.
For any $U\subseteq V(\Gamma)$ with $U\neq\emptyset$, 
$\Gamma\sm U$ has at most one 2-facet, thus
$\pure_2(\Gamma\sm U)$ is shellable.
What remains to do is to determine 
the minimal sets $A$ of edges between $\sigma_1$ and $\sigma_2$
such that 
$\pure_1(\Gamma\sm U)$ is connected for all $U\subseteq V(\Gamma)$ with 
$U\neq \emptyset$.

First, we consider Case 1.
If the set $A$ of edges between $\sigma_1$ and $\sigma_2$ 
shares only one vertex with $\sigma_i$ ($i=1$ or $2$),
then deleting the vertex makes the pure $1$-skeleton disconnected.
Thus $A$ shares at least two vertices with each of $\sigma_1$ and $\sigma_2$.
Without loss of generality we can assume that $A$ contains 
two edges $\{a,d\}$ and $\{b,e\}$.
Consider to take $U=\{a,e\}$.
In order for $\pure_1(\Gamma\sm U)$ to be connected,
one of the edges $\{c,f\}$, $\{c,d\}$ or $\{b,d\}$ should be contained in $A$ (upto isomorphism).
When $\{c,f\}$ is in $A$, then the three edges $\{a,d\}$, $\{b,e\}$ and $\{c,f\}$ 
together with $\sigma_1$ and $\sigma_2$
form (1a) of Figure~\ref{fig:obstructions}, and it is easy to check this is
an edge-minimal obstruction to shellability.
When $\{c,d\}$ or $\{b,d\}$ is in $A$, we then consider to take $U=\{b,d\}$.
In order for $\pure_1(\Gamma\sm U)$ to be connected, at least one more edge 
should be contained in $A$ that connects the edges $\{a,c\}$ and $\{e,f\}$.
In both cases, by checking 4 possible ways of adding this additional edge
(as is depicted in the figure),
we observe that the complex contains (1a), (1b), or (1c) of Figure~\ref{fig:obstructions}.
Here, it is easy to check that (1b) and (1c) are also 
edge-minimal obstructions to shellability.
\begin{figure}
\begin{center}
\input{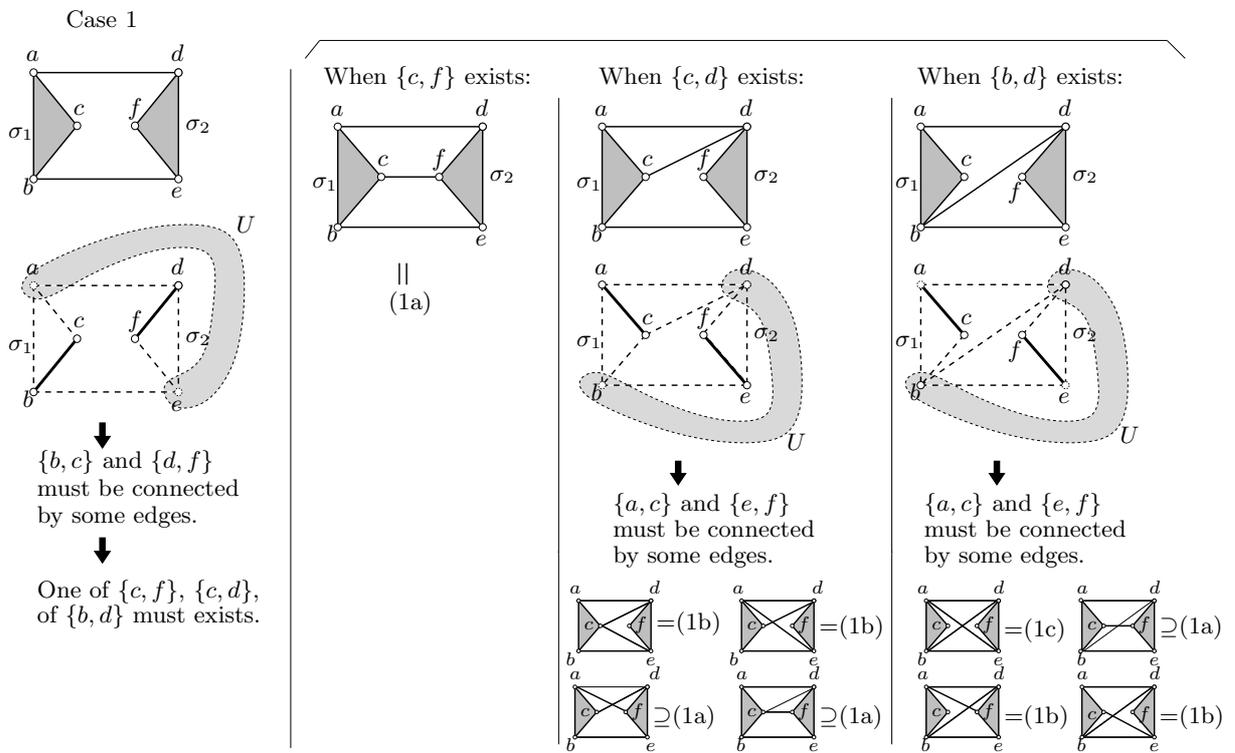}
\end{center}
\caption{Case 1 in the proof of Lemma~\ref{lemma:boundarycase}.}
\label{fig:case1}
\end{figure}

For Case 2, 
consider to take $U=\{a\}$.
In order for $\pure_1(\Gamma\sm U)$ to be connected,
we need an edge between the edges $\{b,c\}$ and $\{d,e\}$.
This gives (2) of Figure~\ref{fig:obstructions}.
It is easy to check this is an edge-minimal obstruction to shellability.
\end{proof}

\smallskip

Next we consider the case there exist nonboundary edges in $\Gamma$.
We divide this case into two subcases by a technical condition. 
The first subcase is Lemma~\ref{lemma:subcase1}, 
and the other is Lemma~\ref{lemma:subcase2}.

\begin{lemma} \label{lemma:subcase1}
Let $\Gamma$ be a $2$-dimensional edge-minimal obstruction to shellability.
Assume there exists a vertex $v$ such that
$\link_\Gamma(\{v\})$ contains a nonboundary edge of $\Gamma$ and
$\pure_1(\link_\Gamma(\{v\}))$ is disconnected.
Then $\Gamma$ is isomorphic to one of (3a), (3b), (3c), (3d) or (3e)
of Figure~\ref{fig:obstructions}.
\end{lemma}
\begin{proof}
Let the nonboundary edge in $\link_\Gamma(\{v\})$ be $\{x,y\}$.
By the assumption that $\pure_1(\link_\Gamma(\{v\}))$ is disconnected,
there is an edge $\{w,u\}$ of $\link_\Gamma(\{v\})$
contained in a different connected component from $\{x,y\}$.
We here have $\Gamma[\{v,x,y,w,u\}]$ nonshellable by Proposition~\ref{prop:link} since
$\pure_1(\link_{\Gamma[\{v,x,y,w,u\}]}(\{v\}))$ is disconnected.
Thus the set of vertices of $\Gamma$ is exactly $\{v,x,y,w,u\}$
because $\Gamma$ is an obstruction to shellability.
Note especially that the 2-facets containing the vertex $v$ are exactly
the two facets $\{v,x,y\}$ and $\{v,w,u\}$ only.

Now it is easily checked that any ways of adding 2-facets
on the four vertices $x,y,w,u$ such that $\{x,y\}$ is a nonboundary edge
give obstructions to shellability
satisfying the assumptions,
and these are isomorphic to 
the five edge-minimal 2-dimensional obstructions (3a)--(3e) 
of Figure~\ref{fig:obstructions}.
\end{proof}

For the next subcase Lemma~\ref{lemma:subcase2},
we prepare three technical lemmas.

\begin{lemma} \label{lemma:linkv}
Let $\Gamma$ be a $2$-dimensional simplicial complex, and
$\pure_1(\link_\Gamma(\{u\}))$ is connected 
for any vertex $u$ with $\link_\Gamma(\{u\})$ containing nonboundary edges of $\Gamma$.
Assume there is a vertex $v$ such that 
$\Gamma\sm\{v\}$ is shellable, 
$\link_\Gamma(\{v\})$ contains nonboundary edges of $\Gamma$,
and these nonboundary edges are connected.
Then $\Gamma$ is shellable.
\end{lemma}
\begin{proof}
By Proposition~\ref{prop:rearrange}, there is a shelling
$\sigma^{(2)}_1,\dotsc,\sigma^{(2)}_{s_2}$,
$\sigma^{(1)}_1,\dotsc,\sigma^{(1)}_{s_1}$,
$\sigma^{(0)}_1,\dotsc,\sigma^{(0)}_{s_0}$
of $\Gamma\sm \{v\}$, where the superscript of each facet indicates
the dimension of the facet.
On the other hand, 
from the assumption that 
the nonboundary edges in $\link_\Gamma(\{v\})$ are connected and 
$\pure_1(\link_\Gamma(\{v\}))$ is connected,
there is a shelling of $\link_\Gamma(\{v\})$,
$\tau^{(1)}_1,\dotsc,\tau^{(1)}_{t_1}$,
$\tau'^{(1)}_1,\dotsc,\tau'^{(1)}_{t'_1}$,
$\tau^{(0)}_1,\dotsc,\tau^{(0)}_{t_0}$,
where $\tau^{(1)}_i$'s are nonboundary edges, and
$\tau'^{(1)}_j$'s are boundary edges of $\Gamma$.
(Note that $t_1\ge 1$ by the assumption that $\link_\Gamma(\{v\})$ contains
nonboundary edges.)
Now the sequence of facets 
$\sigma^{(2)}_1,\dotsc,\sigma^{(2)}_{s_2},
\tau^{(1)}_1\cup\{v\},\dotsc,\tau^{(1)}_{t_1}\cup\{v\},
\tau'^{(1)}_1\cup\{v\},\dotsc,\tau'^{(1)}_{t'_1}\cup\{v\},
\sigma'^{(1)}_1,\dotsc,\sigma'^{(1)}_{s'_1},
\tau^{(0)}_1\cup\{v\},\dotsc,\tau^{(0)}_{t_0}\cup\{v\},
\sigma'^{(0)}_1,\dotsc,\sigma'^{(0)}_{s'_0}$ becomes a shelling of $\Gamma$,
where $\sigma'^{(1)}_1,\dotsc,\sigma'^{(1)}_{s'_1}$ 
and $\sigma'^{(0)}_1,\dotsc,\sigma'^{(0)}_{s'_0}$
are the subsequences
of $\sigma^{(1)}_1,\dotsc,\sigma^{(1)}_{s_1}$ and
$\sigma^{(0)}_1,\dotsc,\sigma^{(1)}_{s_0}$ that list facets of $\Gamma$ only.

To see this sequence is a shelling, it is convenient to use 
the following property related to the \textit{restriction map}.
Let $F_1,F_2,\dotsc,F_t$ be a sequence of facets of a simplicial complex and
define the restriction map as
${\cal R}(F_i) = \{v\in F_i: F_i\sm\{v\}\in \bigcup_{j=1}^{i-1}\overline{F_j}\}$.
Then the sequence  $F_1,F_2,\dotsc,F_t$ is a shelling of the simplicial complex
if and only if ${\cal R}(F_i)\subseteq F_k$ implies $i\le k$ for all $i$ and $k$.
This property follows from Lemma~2.4 and Proposition~2.5 of \cite{B-W:1996}.

For our sequence of facets, for each facet $\eta$ appearing in the sequence,
we check the condition as follows.
First, for $\eta=\sigma^{(2)}_i$, no violation can occur since
$\sigma^{(2)}_1,\dotsc,\sigma^{(2)}_{s_2}$ is the first part of a shelling of $\Gamma$.
No violation can occur either for 
$\eta=\tau^{(1)}_i\cup\{v\}$, $\sigma'^{(1)}_j$, $\tau^{(0)}_k\cup\{v\}$, 
and $\sigma'^{(0)}_l$, 
since violation can occur only when $\dim\eta - \dim{\cal R}(\eta)\ge 2$.
(When $\eta=\tau^{(1)}_1\cup\{v\}$, ${\cal R}(\eta)=\{v\}$ and we have 
$\dim\eta - \dim{\cal R}(\eta)=2$, but for this case obviously ${\cal R}(\eta)=\{v\}$ is
not contained in any of the previous facets.)
What remains is the case $\eta=\tau'^{(1)}_l\cup\{v\}=\{a,b\}\cup\{v\}$.
For this case we have ${\cal R}(\eta)=\{a\}$ (or $\{b\}$) or $\{a,b\}$,
and we need to check the former.
Assume ${\cal R}(\eta)=\{a\}$ and $\{a\}$ is contained in 
some previous facet $\eta'=\{a,x,y\}$. 
Note that $\eta'$ cannot contain $v$ in order for ${\cal R}(\eta)=\{a\}$, 
thus $\eta'$ should be one of $\sigma^{(2)}_i$, and $x$ and $y$ are distinct from $v$.
We here have that $\link_\Gamma(\{a\})$ contains 
the edge $\{v,b\}$ that is a nonboundary edge of $\Gamma$, 
and that $\pure_1(\link_\Gamma(\{a\}))$
is connected by the assumption.
Hence there is a path $v\mbox{-}p\mbox{-}\dotsm\mbox{-}y$
connecting $\{v,b\}$ and $\{x,y\}$.
(See Figure~\ref{fig:proof_linkv}. 
Note that this path, connecting two edges $\{v,b\}$ and $\{x,y\}$, starts from $v$,
not $b$, by the assumption that the edge $\{a,b\}$ is a boundary edge of $\Gamma$.
The last vertex can be $x$ instead of $y$.)
Here we have that the edge $\{a,p\}$ is a nonboundary edge of $\Gamma$ and 
this means that the facet $\{a,p\}\cup\{v\}$ is listed in our sequence
as $\tau^{(1)}_j\cup\{v\}$ before $\eta$.
This contradicts the assumption that ${\cal R}(\eta)=\{a\}$.
Thus we conclude that $\{a\}$ is not contained in any of the previous facets.

\begin{figure}
\begin{center}
\scalebox{.4}{\input{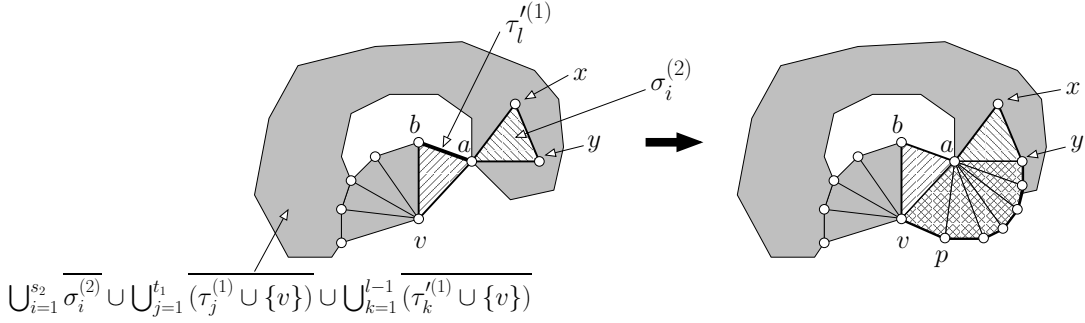}}
\end{center}
\caption{Proof of Lemma~\ref{lemma:linkv}: when a possible violation occurs.}
\label{fig:proof_linkv}
\end{figure}
\end{proof}

\begin{lemma} \label{lemma:disjoint-nbd-link}
Let $\Gamma$ be a $2$-dimensional simplicial complex.
Assume there is a vertex $v$ with $\link_\Gamma(\{v\})$
containing nonboundary edges of $\Gamma$, 
such that $\Gamma\sm\{v\}$ is shellable and 
the nonboundary edges of $\Gamma$ in $\link_\Gamma(\{v\})$ are
not connected.
Then $\Gamma$ is nonshellable.
\end{lemma}
\begin{proof}
Let $\Gamma_1 = \pure_2(\Gamma\sm\{v\})$ and 
$\Gamma_2 = \{v\}*\pure_1(\link_{\Gamma}(\{v\})) 
:= 
\{\{v\}\cup\tau : \tau\in\pure_1(\link_{\Gamma}(\{v\}))\}
\cup
\pure_1(\link_{\Gamma}(\{v\}))$.
(Note that $\Gamma\sm\{v\}$ is 2-dimensional since $\link_\Gamma(\{v\})$
 contains nonboundary edges.)
We have $\Gamma_1\cup\Gamma_2 = \pure_2(\Gamma)$.
From the Mayer-Vietris sequence, we have an exact sequence
$\widetilde{H}_1(\Gamma_1)\oplus\widetilde{H}_1(\Gamma_2) \rightarrow
 \widetilde{H}_1(\Gamma_1\cup\Gamma_2) \rightarrow
 \widetilde{H}_0(\Gamma_1\cap\Gamma_2) \rightarrow
 \widetilde{H}_0(\Gamma_1)\oplus\widetilde{H}_0(\Gamma_2)$.
Since $\Gamma_1=\Gamma\sm\{v\}$ is shellable with $\dim \Gamma_1 = 2$, 
$\widetilde{H}_1(\Gamma_1)=\widetilde{H}_0(\Gamma_1)=0$
by Proposition~\ref{prop:homology}.
Also we have $\widetilde{H}_1(\Gamma_2)=\widetilde{H}_0(\Gamma_2)=0$ 
since $\Gamma_2$ is a cone.
Thus we have $\widetilde{H}_1(\Gamma_1\cup\Gamma_2)
\cong \widetilde{H}_0(\Gamma_1\cap\Gamma_2)$.
On the other hand, the assumption that the nonboundary edges of $\Gamma$ in 
$\link_\Gamma(\{v\})$ are disconnected implies that
$\Gamma_1\cap\Gamma_2$ has at least two connected components.
This implies $\widetilde{H}_0(\Gamma_1\cap\Gamma_2)\neq 0$.
Thus we have 
$\widetilde{H}_1(\pure_2(\Gamma))
=\widetilde{H}_1(\Gamma_1\cup\Gamma_2)
\cong \widetilde{H}_0(\Gamma_1\cap\Gamma_2)\neq 0$.
By Proposition~\ref{prop:homology}, $\pure_2(\Gamma)$ is nonshellable.
Hence $\Gamma$ is nonshellable by Proposition~\ref{prop:purepart}.
\end{proof}

\begin{lemma} \label{lemma:disjoint-twoedges}
Let $\Gamma$ be a $2$-dimensional obstruction to shellability.
Assume that,
for any vertex $v$
with $\link_\Gamma(\{v\})$ containing nonboundary edges,
the nonboundary edges of $\Gamma$ in $\link_\Gamma(\{v\})$ 
are not connected.
Then, for each vertex $v$
with $\link_\Gamma(\{v\})$ containing nonboundary edges,
$\link_\Gamma(\{v\})$ contains
exactly two nonboundary edges disjointly.
\end{lemma}
\begin{proof}
Assume $\link_\Gamma(\{v\})$ contains nonboundary edges.
Since they are not connected, we can take two nonboundary edges
$\{a,b\}$ and $\{c,d\}$ from distinct connected components.
Assume there is another nonboundary edge $\epsilon$ in $\link_\Gamma(\{v\})$.
If, under this assumption, 
we find two vertices, say $v$ and $w$, such that $\link_{\Gamma\sm\{w\}}(\{v\})$
contains two disconnected nonboundary edges,
then this makes a contradiction and completes the proof 
since $\Gamma\sm\{w\}$ becomes nonshellable by Lemma~\ref{lemma:disjoint-nbd-link}
together with the fact $\Gamma\sm\{v,w\}$ is shellable.

First, 
consider the case that the edge $\epsilon$ does not share any vertex with
$\{a,b\}$ and $\{c,d\}$.
Assume $\epsilon=\{e,f\}$, where $e$ and $f$ are distinct from $a,b,c,d$,
and that $\epsilon$ is in a different connected component than
that of $\{a,b\}$.
(We do not assume whether or not
$\epsilon$ is in the component that $\{c,d\}$ belongs to.)
Then, we can pick a vertex $w\in \{c,d,e,f\}$ 
such that $\link_{\Gamma\sm\{w\}}(\{v\})$ contains $\{a,b\}$ and
one of $\{c,d\}$ and $\{e,f\}$ as nonboundary edges of $\Gamma\sm\{w\}$ and we are done.
This can be confirmed as follows.
If a nonboundary edge $\eta$ ($\in \{\{a,b\},\{c,d\},\{e,f\}\}$)
becomes a boundary edge by deleting a vertex $w$
($w\notin\eta$),
then the edge $\eta$ belongs to exactly one 2-facet
$\{w\}\cup \eta$ other than $\{v\}\cup \eta$.
Since we have four candidates $c,d,e,f$ for $w$
while there are only three edges $\{a,b\}$, $\{c,d\}$, $\{e,f\}$
for the edge $\eta$, we can always pick a free vertex $w$ from $\{c,d,e,f\}$
such that deleting $w$ does not affect these nonboundary edges
except the one $w$ belongs to.

Next, consider the case that the edge $\epsilon$ shares a vertex with
$\{a,b\}$ or $\{c,d\}$. Without loss of generality we assume $\epsilon=\{c,e\}$.
Note that $e$ is distinct from $a$ and $b$.
If the $2$-facet $\{c,d,e\}$ does not exist,
then deleting $d$ or $e$
preserves the edge $\{c,e\}$ or $\{c,d\}$ as a nonboundary edge, respectively.
Further, deleting one of $d$ or $e$ preserves the edge $\{a,b\}$
as a nonboundary edge by a similar argument as above.
Hence we can take a vertex $w\in\{d,e\}$ such that
$\link_{\Gamma\sm\{w\}}(\{v\})$ contains two nonboundary edges
disjointly and we are done.
Thus we assume that the $2$-facet $\{c,d,e\}$ exists.
Since $\{v,a,b,c,d\}$ is a proper subset of $V(\Gamma)$,
$\Gamma[\{v,a,b,c,d\}]$ is shellable,
and thus $\pure_1(\link_{\Gamma[\{v,a,b,c,d\}]}(\{v\}))$ is connected
by Propositions~\ref{prop:link} and \ref{prop:1dim}.
Without loss of generality, we can assume there is a boundary edge
$\{a,d\}$ or $\{a,c\}$ in $\link_\Gamma(\{v\})$.
Assume the edge $\{a,d\}$ exists in $\link_\Gamma(\{v\})$.
Then in $\link_\Gamma(\{d\})$, there are nonboundary edges
$\{v,a\}$, $\{v,c\}$ and $\{c,e\}$ connected in a path.
Since the nonboundary edges in $\link_\Gamma(\{d\})$ are not connected,
there exists another nonboundary edge $\eta$ in other connected component.
Here we have two nonboundary edges $\eta$ and $\{v,a\}$ in a different
connected component, and there is the third nonboundary edge $\{c,e\}$
not sharing any vertices with $\eta$ and $\{v,a\}$.
Thus we are lead to a contradiction via the first case.
Hence we assume an edge $\{a,c\}$ is in $\link_\Gamma(\{d\})$
instead of $\{a,d\}$.
We have two connected nonboundary edges $\{v,c\}$ and $\{c,e\}$
in $\link_\Gamma(\{d\})$. Since the nonboundary edges in $\link_\Gamma(\{d\})$
are not connected, 
there exists a nonboundary edge $\{f,g\}$
in $\link_\Gamma(\{d\})$, in a different connected component from
$\{v,c\}$ and $\{c,e\}$.
We here have $f$ and $g$ distinct from $v$, $c$, $d$, or $e$.
Since $\{f,g\}$ is a nonboundary edge, there exists a 2-facet $\{f,g,h\}$ in $\Gamma$.
Here, $\link_\Gamma(\{d\})$ contains two disjoint nonboundary edges 
$\{v,c\}$ and $\{f,g\}$, and any deletion of one vertex other than $c$, $d$ and $v$
preserves $\{v,c\}$ as a nonboundary edge since it is contained in
three 2-facets $\{v,a,c\}$, $\{v,c,d\}$ and $\{v,c,e\}$.
Now, if there is a vertex $w$ other than $c$, $d$, $f$, $g$, $h$, and $v$, then
deleting $w$ leaves $\{f,g\}$ and $\{v,c\}$ as two disjoint nonboundary edges 
in $\link_\Gamma(\{d\})$.
Otherwise, we should have $\{f,g\}=\{a,b\}$, and in this case 
deleting the vertex $e$ preserves
the two disjoint edges $\{f,g\}$ and $\{v,c\}$ as nonboundary edges.
\end{proof}

Now we give the last subcase of Theorem~\ref{thm:ots} in the following lemma.

\begin{lemma} \label{lemma:subcase2}
Let $\Gamma$ be an edge-minimal $2$-dimensional obstruction to shellability.
Assume that 
there is a nonboundary edge in $\Gamma$, and that
$\pure_1(\link_\Gamma(\{v\}))$ is connected 
for any vertex $v$ with
$\link_\Gamma(\{v\})$ containing nonboundary edges.
Then $\Gamma$ is isomorphic to one of (4a), (4b) or (4c)
of Figure~\ref{fig:obstructions}.
\end{lemma}

\begin{proof}
First observe that, for every vertex $v$ containing nonboundary edges
of $\Gamma$,
the nonboundary edges in $\link_\Gamma(\{v\})$ are not connected,
since otherwise $\Gamma$ becomes shellable by Lemma~\ref{lemma:linkv}.
(Note that such a vertex $v$ surely exists by the assumption that $\Gamma$ has a nonboundary edge.)
By Lemma~\ref{lemma:disjoint-twoedges},
for such a vertex $v$, $\link_\Gamma(\{v\})$ contains exactly two
nonboundary edges $\{x_1,x_2\}$ and $\{x_3,x_4\}$ disjointly.
Here, if $V(\Gamma) = \{v,x_1,x_2,x_3,x_4\}$, then 
the two nonboundary edges $\{x_1,x_2\}$ and $\{x_3,x_4\}$ are connected
by a boundary edge in $\link_\Gamma(\{v\})$ since
$\pure_1(\link_\Gamma(\{v\}))$ is connected.
On the other hand, if $V(\Gamma) \supsetneq \{v,x_1,x_2,x_3,x_4\}$, then
$\Gamma[\{v,x_1,x_2,x_3,x_4\}]$ is shellable and,
by Propositions~\ref{prop:link} and \ref{prop:1dim},
$\pure_1(\link_{\Gamma[\{v,x_1,x_2,x_3,x_4\}]}(\{v\}))$ is connected.
Thus again the two nonboundary edges $\{x_1,x_2\}$ and $\{x_3,x_4\}$ 
are connected by a boundary edge in $\link_\Gamma(\{v\})$.
Hence we have that the two nonboundary edges 
$\{x_1,x_2\}$ and $\{x_3,x_4\}$ are always connected
by a boundary edge in $\link_\Gamma(\{v\})$.

We have another consequence.
In $\link_\Gamma(\{v\})$, we have two nonboundary edges $\{x_1,x_2\}$ and $\{x_3,x_4\}$ disjointly,
thus there are 2-facets $\{x_1,x_2,x_5\}$ and $\{x_3,x_4,x_6\}$.
Here, $v,x_1,x_2,x_3,x_4$ are all distinct, but $x_5$ and $x_6$  can be 
identical to some of $x_1,x_2,x_3,x_4$.
Since $\Gamma$ is an obstruction to shellability and
$\{x_1,x_2,x_3,x_4,x_5,x_6\}$ is a proper subset of $V(\Gamma)$,
$\Gamma[\{x_1,x_2,x_3,x_4,x_5,x_6\}]$ is shellable.
This fact together with the fact that the nonboundary edges 
$\{x_1,x_2\}$ and $\{x_3,x_4\}$ of $\Gamma$ in $\link_\Gamma(v)$ are disconnected
implies that $\Gamma[\{v,x_1,x_2,x_3,x_4,x_5,x_6\}]$ is nonshellable
by Lemma~\ref{lemma:disjoint-nbd-link}.
Hence we have $V(\Gamma) = \{v,x_1,x_2,x_3,x_4,x_5,x_6\}$.
Thus we have $\lvert V(\Gamma)\rvert = 5,6$ or $7$.
(This together with Lemmas~\ref{lemma:boundarycase} and \ref{lemma:subcase1}
 gives a different proof of Theorem~\ref{thm:wachs} (iv).)

\medskip
Now pick a vertex $v_1$ of $\Gamma$, and
assume that $\{x,y\}$ is a nonboundary edge of $\Gamma$ in 
$\link_\Gamma(\{v_1\})$.
Then there exists a boundary edge adjacent to $\{x,y\}$ in 
$\link_\Gamma(\{v_1\})$.
We assume this boundary edge is $\{x,v_2\}$ without loss of generality. 
Here the edge $\{x,v_1\}$ is a nonboundary edge in $\link_\Gamma(\{v_2\})$. 
This requires the existence of the boundary edge $\{v_1,v_3\}$ in 
$\link_\Gamma(\{v_2\})$.
(Remark that the boundary edge cannot be $\{x,v_3\}$ instead of $\{v_1,v_3\}$
because $\{x,v_2\}$ is a boundary edge.)
\begin{figure}
\begin{center}
\input{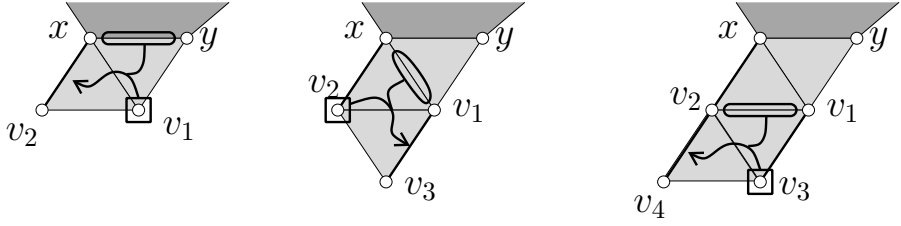}
\end{center}
\caption{Proof of Lemma~\ref{lemma:subcase2}: a band grows from $x v_1$.}
\label{fig:band}
\end{figure}
Repeating this procedure makes a band grown from the edge $x v_1$
as indicated in Figure~\ref{fig:band}.
But, since the number of vertices of $\Gamma$ is at most 7, 
there is a number $N$ with $v_N = v_i$ for some $i<N$.
By renaming the vertices, this ensures that 
there is a subcomplex $\Delta'$ on vertices $v_1,v_2,\dotsc,v_{n}$ 
(all distinct) whose facets are
$\{\{v_1,v_2,v_3\},\{v_2,v_3,v_4\},\dotsc,\{v_{n-1},v_n,v_1\}\}$,
where the edges $\{v_1,v_3\},\{v_2,v_4\},\dotsc,\{v_{n-2},v_n\}$ and
$\{v_{n-1},v_1\}$ are boundary edges of $\Gamma$.
(See the left of Figure~\ref{fig:preband}.)
We have $n\ge 5$ in order for the required edges to be boundary edges.
\begin{figure}
\begin{center}
\input{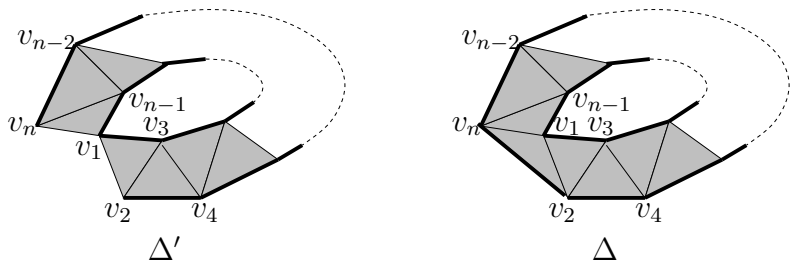}
\caption{The subcomplexes $\Delta'$ and $\Delta$. (The bold edges are boundary edges
 of $\Gamma$.)}
\label{fig:preband}
\end{center}
\end{figure}

Here, since the vertex $v_1$ has two disjoint nonboundary edges
$\{v_2,v_3\}$ and $\{v_{n-1},v_n\}$ in $\link_\Gamma(\{v_1\})$,
these two nonboundary edges must be connected by a boundary edge.
Only possible choice for this boundary edge is $\{v_n,v_2\}$,
since $\{v_{n-1},v_1\}$ and $\{v_1,v_3\}$ are boundary edges.
Thus there is a subcomplex $\Delta$
on vertices $\{v_1,v_2,\dotsc,v_n\}$ whose facets are
$\{v_1,v_2,v_3\},\{v_2,v_3,v_4\},\dotsc,\{v_{n-1},v_n,v_1\}$,
$\{v_n,v_1,v_2\}$,
where the edges $\{v_1,v_3\}$, $\{v_2,v_4\},\dotsc,\{v_{n-2},v_n\}$,
$\{v_{n-1},v_1\}, \{v_n,v_2\}$ are boundary edges of $\Gamma$.
This $\Delta$ is a cylinder or a M\"obius band according to the parity of $n$,
see the right of Figure~\ref{fig:preband}.
Since we have $5\le n\le 7$, $\Delta$ forms one of the simplicial complexes
of (4a), (4b) or (4c) in Figure~\ref{fig:obstructions}.

Assume there is an extra 2-facet 
in $\Gamma[\{v_1,v_2,\dotsc,v_n\}]$.
Since each $v_i$ already has two nonboundary edges in 
$\link_{\Gamma[\{v_1,v_2,\dotsc,v_n\}]}(\{v_i\})$,
the extra $2$-facet can not have nonboundary edges
by Lemma~\ref{lemma:disjoint-twoedges}.
Thus all the three edges of the extra facet are boundary edges.
But this is impossible, since no three nonedges form
a triangle in the three complexes (4a)-(4c) of Figure~\ref{fig:obstructions}.
Hence we 
exactly
have $\pure_2(\Gamma[\{v_1,v_2,\dotsc,v_n\}])=\Delta$.
Here we have $\Delta$ nonshellable by Proposition~\ref{prop:homology}
because $\widetilde{H}_1(\Delta)\cong\mathbb{Z}$,
thus $\Gamma$ can not have extra vertices than $\{v_1,v_2,\dotsc,v_n\}$
because $\Gamma$ is an obstruction to shellability.

Here, already in \cite[Theorem 5]{W:2000} 
it is pointed out that
these three complexes (4a)-(4c) of Figure~\ref{fig:obstructions}
themselves are obstructions to shellability.
Hence we conclude that 
the 2-dimensional edge-minimal obstruction to shellability $\Gamma$
is one of the three complexes.

\end{proof}

\begin{proof}[Proof of Theorem~\ref{thm:ots}]
Lemmas~\ref{lemma:boundarycase}, \ref{lemma:subcase1}
and \ref{lemma:subcase2} show the statement.
\end{proof}

\section{Obstructions to partitionability and sequential Cohen-Macaulayness}
\label{sec:otp-otCM}

Partitionability and sequential Cohen-Macaulayness are well-known properties
implied by shellability. 

\begin{defn}
A simplicial complex $\Gamma$ is \textit{partitionable} if $\Gamma$ can be partitioned
into intervals as
$\Gamma = \dot\bigcup[\tau_\sigma, \sigma],$
where 
$[\tau_\sigma, \sigma] = \{\eta\in\Gamma: \tau_\sigma\subseteq\eta\subseteq\sigma\}$
is the interval between a facet $\sigma$ and one of its face $\tau_\sigma$.
\end{defn}

\begin{defn}\label{defn:sCM}
\begin{enumerate}
\setlength{\itemsep}{0pt}
\item
A pure simplicial complex $\Gamma$ is \textit{Cohen-Macaulay} (over $\mathbb{Z}$) if,
for any $\tau\in\Gamma$,
$\widetilde{H}_k(\link_\Gamma(\tau)) = 0$ unless $k=\dim\link_\Gamma(\tau)$,
where $\widetilde{H}_k$ is the $k$-th reduced homology group
(over $\mathbb{Z}$).
\item
A simplicial complex $\Gamma$ is \textit{sequentially Cohen-Macaulay} (over $\mathbb{Z}$)
if $\pure_i(\Gamma)$ is Cohen-Macaulay (over $\mathbb{Z}$) 
for all $0\le i\le \dim\Gamma$.
\end{enumerate}
\end{defn}

See \cite{B-W:1996}, \cite{KO:1996}, etc.\ for partitionability.
Sequential Cohen-Macaulayness is originally defined in terms of
commutative algebra by Stanley, see \cite[Section~III.2]{St:1996}.
In Definition~\ref{defn:sCM} above, (i) is a characterization given by
Reisner~\cite{R:1976}, and (ii) is a characterization given by
Duval~\cite{D:1996}.
See also \cite{BWV:2009} for sequential Cohen-Macaulayness.
Basic properties about partitionability and sequential Cohen-Macaulayness
we need in this paper are as follows.

\begin{prop}
\label{prop:implication}
\begin{enumerate}
\setlength{\itemsep}{0pt}
\item If $\Gamma$ is shellable, then $\Gamma$ is partitionable.
\item If $\Gamma$ is shellable, then $\Gamma$ is sequentially Cohen-Macaulay.
\end{enumerate}
\end{prop}
\begin{proof}
(i) is shown in \cite[Section~2]{B-W:1996}. (ii) appears in 
\cite{BWV:2009}, \cite{D:1996}, \cite[Section~III.2]{St:1996}, \cite{W:1999}, etc.
\end{proof}

(Remark: It is known that partitionability does not imply 
sequential Cohen-Macaulayness. 
On the other hand, it is an open problem whether sequential Cohen-Macaulayness 
implies partitionability or not. See \cite{D:1996}, \cite[Section~III.2]{St:1996}, etc.)

\begin{prop}\label{prop:link-psCM}
\begin{enumerate}
\setlength{\itemsep}{0pt}
\item 
If $\Gamma$ is partitionable, then $\link_\Gamma(\tau)$ is partitionable
for any $\tau\in\Gamma$.
\item 
If $\Gamma$ is sequentially Cohen-Macaulay, 
then $\link_\Gamma(\tau)$ is sequentially Cohen-Macaulay
for any $\tau\in\Gamma$.
\end{enumerate}
\end{prop}
\begin{proof}
(i) is essentially given by \cite[Proposition~3.1]{KO:1996}. Though the discussion
given there is restricted only for pure cases, but it naturally applies for
nonpure cases. 
(ii) follows immediately from the equivalent definition given 
in \cite[Definition~1.2]{BWV:2009}.
\end{proof}

\begin{prop} \label{prop:01dim-psCM}
\begin{enumerate}
\setlength{\itemsep}{0pt}
\item
$0$-dimensional simplicial complexes are always 
partitionable and sequentially Cohen-Macaulay.
\item
A $1$-dimensional simplicial complex $\Gamma$ is sequentially Cohen-Macaulay
if and only if $\pure_1(\Gamma)$ is connected.
\item
A $1$-dimensional simplicial complex $\Gamma$ is partitionable if and only if
$\pure_1(\Gamma)$ has at most one connected component that is a tree 
(i.e., component with no cycle).
\end{enumerate}
\end{prop}
\begin{proof}
(i) and (ii) are clear from the definition.
(iii) is essentially pointed out in \cite[Section~36]{KP:2003}.
\end{proof}

\begin{prop} (\cite[Proposition~11.7]{B:1995}) \label{prop:strconn}
If a pure $d$-dimensional simplicial complex is Cohen-Macaulay, then
it is strongly connected, i.e.,
for any two facets $\sigma$ and $\sigma'$, 
there is a sequence of facets $\sigma=\sigma_1,\sigma_2,\dotsc,\sigma_k=\sigma'$
such that $\sigma_i\cap\sigma_{i+1}$ is a $(d-1)$-dimensional face of $\Gamma$
for all $1\le i\le k-1$.
\end{prop}

\begin{prop} \label{prop:2facets-nonp}
If  a simplicial complex $\Gamma$ has two facets $\sigma_1,\sigma_2$ 
with $\dim\sigma_1,\dim\sigma_2\ge 1$
such that each of the subfaces of $\sigma_i$ of dimension $\dim\sigma_i-1$
belongs only to $\sigma_i$ ($i=1,2$), then $\Gamma$ is not partitionable.
\end{prop}
\begin{proof}
For such a facet $\sigma_i$ as in the statement, $\tau_{\sigma_i}$ 
should be $\emptyset$
in order that all its subfaces of dimension $\dim\sigma_i-1$ are covered 
in the partition.
Thus there exists no partition if there are more than one such facets.
\end{proof}

\begin{prop} \label{prop:band-nonp}
Let $\Gamma$ be a $d$-dimensional simplicial complex on vertices
$\{v_0,v_1,\dotsc, v_{n-1}\}$ with $d\ge 2$.
Assume $d$-facets 
$\{\{v_0,v_1,\dotsc,v_{d}\}, \{v_1,v_2,\dotsc,v_{d+1}\},\dotsc,
\{v_{n-1},v_0,\dotsc,v_{d-1}\}\}$ are contained in $\Gamma$, and
among the $(d-1)$-dimensional subfaces of these $d$-facets,
each of $\{v_0,v_1,\dotsc,v_{d-1}\}, \{v_1,v_2,\dotsc,v_{d}\},\dotsc,
\{v_{n-1},v_0,\dotsc,v_{d-2}\}$ belongs to exactly two $d$-facets
and each of other $(d-1)$-subfaces belongs to 
only one $d$-facet. 
(For this, we need $n\ge 2d+1$.)
Then $\Gamma$ is not partitionable.
\end{prop}
\begin{proof}
Throughout this proof, addition in the index is considered as modulo $n$.
We denote $\sigma_k = \{v_k,v_{k+1},\dotsc,v_{k+d}\}$.

Assume $\Gamma$ is partitionable and has a partition
$\Gamma = \dot\bigcup [\tau_\sigma, \sigma]$, where $\sigma$ is a facet of $\Gamma$
and $\tau_\sigma$ is a subface of $\sigma$.
Let us denote $\tau_{\sigma_k}=\tau_k$ for convenience.
Since
the $(d-1)$-faces
$\sigma_k \sm \{v_{k+j}\}$
($1\le j\le d-1$)
belong only to
the $d$-facet $\sigma_k$,
$\tau_k$ should be one of
$\emptyset$, $\{v_k\}$, $\{v_{k+d}\}$, or $\{v_k,v_{k+d}\}$
for each $1\le k\le n$.

Here we observe that $v_k\in \tau_k$ implies $v_{k+d+1}\not\in \tau_{k+1}$,
since otherwise the face $\{v_{k+1},\dotsc,v_{k+d}\}$ is not contained in 
both of 
$[\tau_k,\sigma_k]$ and $[\tau_{k+1},\sigma_{k+1}]$ and thus not contained in
any of the intervals.
Symmetrically, we also have that $v_{k+d}\in \tau_k$ implies $v_{k-1}\not\in\tau_{k-1}$.
Now assume that $\tau_k=\{v_k,v_{k+d}\}$. Then, from these two observations,
we have that $\tau_{k+1}\subseteq\{v_{k+1}\}$ and 
that $\tau_{k-1}\subseteq\{v_{k+d-1}\}$.
This implies that $\{v_{k+1},v_{k+d-1}\}$ is contained 
in both of $[\tau_{k+1},\sigma_{k+1}]$
and $[\tau_{k-1},\sigma_{k-1}]$, a contradiction.
Hence we have that $\tau_k$ is one of
$\emptyset$, $\{v_k\}$, or $\{v_{k+d}\}$, for each $1\le k\le n$.

Assume $\tau_k=\emptyset$ for some $k$.
We have that $\tau_i$ is a single vertex for $i=0,1,\dotsc,n-1$ with $i\neq k$, and
these vertices are distinct.
But $[\tau_k,\sigma_k]$ covers $d+1$ vertices and thus only $n-d-1$ vertices are left,
a contradiction. 

Now we have that $\tau_k$ is one of $\{v_k\}$ or $\{v_{k+d}\}$, for each $1\le k\le n$.
This implies that each interval $[\tau_i,\sigma_i]$ contains exactly one distinct vertex,
thus all the $n$ vertices $\{v_0\},\{v_1\},\dotsc,\{v_{n-1}\}$ are covered by
the $n$ intervals $[\tau_0,\sigma_0], \dotsc, [\tau_{n-1},\sigma_{n-1}]$.
Since any of the intervals $[\tau_k,\sigma_k]$ contains $\emptyset$, 
there should exist some facet $\sigma$ 
other than the $n$ facets $\sigma_0,\sigma_1,\dotsc,\sigma_{n-1}$
such that $\tau_\sigma=\emptyset$.
But, this is impossible since
$\Gamma$ has only $n$ vertices and
they are already covered by the intervals. 
\end{proof}

\bigskip
In this section we discuss the relation between the sets of obstructions 
to the three properties,
shellability, partitionability and sequential Cohen-Macaulayness.
The key of our discussion is the following proposition.

\begin{prop} \label{prop:difference}
Let $\cal P$ and $\cal Q$ be properties of simplicial complexes
such that $\cal P$ implies $\cal Q$.
In a class $\cal X$ of simplicial complexes that is closed under
restriction,
the following are equivalent.
\begin{enumerate}
\setlength{\itemsep}{0pt}
\item
In the class $\cal X$, the set of obstructions to $\cal P$
is different from the set of obstructions to $\cal Q$.
\item
There exists an obstruction to $\cal P$ in $\cal X$
that is not an obstruction to $\cal Q$.
\item
There exists an obstruction to $\cal P$ in $\cal X$
that satisfies $\cal Q$.
\end{enumerate}
\end{prop}
\begin{proof}
\noindent
The implications (iii)$\Rightarrow$(ii)$\Rightarrow$(i) are trivial.
We show the implication (i)$\Rightarrow$(iii).

If there is an obstruction to $\cal P$ in $\cal X$
that is not an obstruction to $\cal Q$,
then it does not satisfy $\cal Q$ and we are done.
Assume $\Gamma$ is an obstruction to $\cal Q$ in $\cal X$
that is not an obstruction to $\cal P$.
Since $\Gamma$ does not satisfy $\cal Q$, it does not satisfy $\cal P$, either.
In order that $\Gamma$ is not an obstruction to $\cal P$, 
there exists some 
$W\subsetneq V(\Gamma)$ such that $\Gamma[W]$ does not satisfy $\cal P$.
If we choose such $W$ minimal with respect to inclusion, 
such $\Gamma[W]$ is an obstruction to $\cal P$ that satisfies $\cal Q$.
Note that $\Gamma[W]$ is in $\cal X$ since $\cal X$ is closed under restriction.
\end{proof}

By this proposition
together with Theorem~\ref{thm:ots}
we have the following theorem.

\begin{thm} \label{thm:otp-otCM}
For the class of simplicial complexes of dimension $\le 2$, 
the set of obstructions to partitionability,
that to sequential Cohen-Macaulayness,
and that to shellability all coincide.
\end{thm}

\begin{proof}
We check that each obstruction to shellability
given by Theorem~\ref{thm:ots} is neither partitionable
nor sequentially Cohen-Macaulay.
Then 
Proposition~\ref{prop:difference}
together with 
Proposition~\ref{prop:implication}
shows the statement,
since the class of simplicial complexes of dimension $\le 2$ is closed
under restriction.

The one dimensional obstruction to shellability $2K_2$
is neither partitionable nor sequentially Cohen-Macaulay by
Proposition~\ref{prop:01dim-psCM}.
For the obstructions derived from (1a)-(1c) of Figure~\ref{fig:obstructions},
they are not partitionable by Proposition~\ref{prop:2facets-nonp}.
These are not sequentially Cohen-Macaulay, either,
because their pure $2$-skeletons are not Cohen-Macaulay
by Proposition~\ref{prop:strconn}.
For the obstructions derived from (2) and (3a)-(3e) of Figure~\ref{fig:obstructions},
there is a vertex such that its link has two disjoint edges.
Hence
these obstructions are neither partitionable nor sequentially Cohen-Macaulay
by Propositions~\ref{prop:link-psCM} and \ref{prop:01dim-psCM}.
For the obstructions derived from (4a)-(4c) of Figure~\ref{fig:obstructions},
they are not partitionable by Proposition~\ref{prop:band-nonp},
and not sequentially Cohen-Macaulay since their pure $2$-skeletons have
$\widetilde{H}_1 \cong \mathbb{Z}$.
\end{proof}

\begin{remark}
Both of partitionability and sequential Cohen-Macaulayness
are known to be strictly weaker than shellability already in dimension $2$.
There are nonshellable but partitionable complexes for dimensions 
$\ge 1$,
since the pure 1-skeleton of a 1-dimensional partitionable complex
can be disconnected
by Proposition~\ref{prop:01dim-psCM}
while 
the pure 1-skeleton of a shellable 1-dimensional complex should be connected
by Proposition~\ref{prop:1dim}.
Also there are nonshellable but sequentially Cohen-Macaulay complexes
for dimensions $\ge 2$, see \cite[Section~III.2]{St:1996}.
(Cones over the examples give counterexamples of higher dimensions 
for both.)
But the theorem above does not lead to a contradiction
because these three properties are not closed under restriction.
\end{remark}

\medskip
On the other hand, Woodroofe~\cite{W:2009} showed the following result
that determines obstructions to shellability in the class of
flag complexes.

\begin{thm}(Woodroofe~\cite{W:2009}) \label{thm:woodroofe}
In the class of flag complexes, the obstructions to shellability are
exactly the independence complexes of $C_n$ with $n=4$ or $n\ge 6$,
where $C_n$ is a cycle graph of order $n$.
\end{thm}
Here, an \textit{independence complex} is a simplicial complex
derived from a graph $G$, whose vertices are the vertices of $G$ and
faces are independent sets of $G$.
The independence complex of $C_n$ is a 
$d$-dimensional simplicial complex
with $d=\lfloor\frac{n}{2}\rfloor -1$.
When $n$ is even, 
it has two $d$-facets disjointly.
When $n$ is odd, 
its $d$-facets are 
$\{v_0,v_1,\dotsc,v_{d}\},\{v_1,v_2,\dotsc,v_{d+1}\},\dotsc,
\{v_{n-1},v_0,v_1,\dotsc,v_{d-1}\}$ by suitably indexing the vertices.

Since the class of flag complexes is closed under restriction, 
Proposition~\ref{prop:implication} also gives the following.

\begin{thm} \label{thm:flag}
In the class of flag complexes, the set of obstructions to partitionability,
that to sequential Cohen-Macaulayness, and that to shellability all coincide.
\end{thm}
\begin{proof}
We check the obstructions to shellability given by Theorem~\ref{thm:woodroofe}
are neither sequentially Cohen-Macaulay nor partitionable.
Then we have the statement from Proposition~\ref{prop:implication} 
by setting ${\cal X}$ to be the class of flag complexes,
since the class of flag complexes is closed under restriction.

When $n$ is even with $n\ge 4$, the independence complexes $\Gamma(C_n)$ of $C_n$ 
is not sequentially Cohen-Macaulay because $\pure_d(\Gamma(C_n))$
is not Cohen-Macaulay by Proposition~\ref{prop:strconn},
and it is nonpartitionable by Proposition~\ref{prop:2facets-nonp}.
When $n$ is odd with $n\ge 7$, $\Gamma(C_n)$ is not sequentially Cohen-Macaulay
since $\widetilde{H}_1(\pure_d(\Gamma(C_n))) \cong \mathbb{Z}$,
and its nonpartitionability is shown by Proposition~\ref{prop:band-nonp}.
\end{proof}

\medskip
One application of the discussions on obstructions to
a property $\cal P$ is the following relation to the concept
of ``hereditary-$\cal P$'' property.

\begin{defn} \label{def:hereditary}
Let $\cal P$ be a property of simplicial complexes.
We say that 
a simplicial complex $\Gamma$ is \textit{hereditary-$\cal P$}
if all the restrictions of $\Gamma$ (including $\Gamma$ itself)
satisfy $\cal P$.
\end{defn}

This hereditary-$\cal P$ property can be characterized by
the obstructions to $\cal P$ as follows.

\begin{prop} \label{prop:hereditary}
A simplicial complex $\Gamma$ is \textit{hereditary-$\cal P$}
if and only if $\Gamma[W]$ is not an obstruction to $\cal P$
for all $W\subseteq V(\Gamma)$.
\end{prop}

The proof of this proposition is straight-forward.
By this proposition we can view the obstructions to $\cal P$
as the excluded minors of hereditary-$\cal P$ property
with respect to restriction.
We have the following corollary from Theorems~\ref{thm:otp-otCM} and 
\ref{thm:flag}.

\begin{cor} \label{cor:hereditary}
The following are equivalent
for a simplicial complex $\Gamma$, 
if the dimension of $\Gamma$ is at most $2$, or $\Gamma$ is a flag complex.
\begin{enumerate}
\setlength{\itemsep}{0pt}
\item $\Gamma$ is hereditary-shellable.
\item $\Gamma$ is hereditary-partitionable.
\item $\Gamma$ is hereditary-sequentially Cohen-Macaulay.
\end{enumerate}
\end{cor}

\bigskip
\begin{remark}
There are more properties related to shellability 
with the following implications.

\smallskip
\hspace*{2cm}
vertex decomposable $\Rightarrow$
shellable $\Rightarrow$
constructible $\Rightarrow$ 
\\
\hspace*{5cm}
$\Rightarrow$ sequentially homotopy-CM
$\Rightarrow$
sequentially Cohen-Macaulay
\smallskip

The same implications for pure complexes (vertex decomposable $\Rightarrow$
shellable $\Rightarrow$
constructible $\Rightarrow$ homotopy-CM $\Rightarrow$ Cohen-Macaulay)
are classical and well-known, see \cite[Sections~11.2 and 11.5]{B:1995}.
Recently, Jonsson~\cite[Section~3.6]{J:2008} gave 
a generalized definition of constructibility
that also applies for nonpure complexes,
and completed the implications of general simplicial complexes shown above.
(He calls the generalized definition ``semipure constructibility'' in order to
discriminate from the classical definition that applies only for pure cases.)
For properties of simplicial complexes satisfying 
${\cal P}\Rightarrow {\cal Q}\Rightarrow {\cal R}$,
it is easy to verify 
using Proposition~\ref{prop:difference}
that if the set of obstructions to $\cal P$
is the same as the set of obstructions to $\cal R$
then also the set of obstructions to $\cal Q$ is the same.
Hence the set of obstructions to constructibility and
that to sequential homotopy-CMness coincide with that to shellability
in the class of simplicial complexes of dimensions $\le 2$
and in the class of flag complexes.
Accordingly, hereditary-constructibility and hereditary-sequential
homotopy-CMness are equivalent to hereditary-shellability
in these classes.

On the other hand, the situation is different for vertex decomposability.
In \cite{HK:2009}, a $2$-dimensional simplicial complex is presented that is
shellable but not vertex decomposable,
and any of whose restrictions are shellable.
This implies the existence of a $2$-dimensional
obstruction to vertex decomposability that is shellable.
Hence the set of obstructions to vertex decomposability and
the set of obstructions to shellability differ in dimension $2$.
On the other hand, the proof of Theorem~\ref{thm:woodroofe} given in \cite{W:2009}
implies that the set of obstructions to vertex decomposability is the same as
the set of obstructions to shellability in the class of flag complexes.
\end{remark}

\section{Strong obstructions: toward higher dimensions}
\label{sec:s-obstruction}

In this section we introduce the following definition.

\begin{defn} \label{def:s-obstruction}
Let $\cal P$ be a property of simplicial complexes.
A simplicial complex $\Gamma$ is a \textit{strong obstruction} to ${\cal P}$
if $\Gamma$ does not satisfy $\cal P$
but 
$\link_{\Gamma[W]}(\tau)$ satisfies $\cal P$ for any $W\subseteq V(\Gamma)$ and
$\tau\in\Gamma$ unless $W=V(\Gamma)$ and $\tau=\emptyset$
(i.e., unless $\link_{\Gamma[W]}(\tau)\neq \Gamma$).
\end{defn}

We say a property $\cal P$ is \textit{link-preserving} if
$\link_\Gamma(\tau)$ satisfies $\cal P$ for all $\tau\in\Gamma$  
whenever $\Gamma$ satisfies $\cal P$.
One of the important properties of shellability is that
it is link-preserving as shown in Proposition~\ref{prop:link}.
Partitionability and sequential Cohen-Macaulayness also have
this property, see Proposition~\ref{prop:link-psCM}.
For a link-preserving property $\cal P$, a slightly weaker condition
characterizes strong obstructions as follows.
The proof is straightforward.
\begin{prop}
If $\cal P$ is a link-preserving property,
then a simplicial complex $\Gamma$ is a strong obstruction to $\cal P$
if and only if $\Gamma$ does not satisfy $\cal P$,
$\Gamma[W]$ satisfies $\cal P$ for any $W\subsetneq V(\Gamma)$,
and $\link_\Gamma(\tau)$ satisfies $\cal P$ for any $\tau\in\Gamma$
with $\tau\neq\emptyset$.
\end{prop}

We have the following characterization of hereditary-$\cal P$ property
instead of Proposition~\ref{prop:hereditary}
when the property $\cal P$ is link-preserving.
The proof is immediate using the property that 
the operations of restricting and taking link commute.

\begin{prop} \label{prop:s-hereditary}
Let $\cal P$ be a property of simplicial complexes that is link-preserving.
A simplicial complex $\Gamma$ is hereditary-$\cal P$
if and only if 
$\link_{\Gamma[W]}(\tau)$ is not a strong obstructions to $\cal P$ 
for any $W\subseteq V(\Gamma)$ and $\tau\in\Gamma[W]$.
\end{prop}

The obstructions to $\cal P$ are the excluded minors
of hereditary-$\cal P$ property with respect to restriction operation,
while the strong obstructions are the excluded minors of 
hereditary-$\cal P$ property
with respect to restriction and link operation
for a link-preserving property $\cal P$.

\medskip
As same as Theorem~\ref{thm:wachs}~(iii), there exist strong obstructions
to shellability in each dimension $d\ge 1$.
The obstructions in the class of flag complexes given by Theorem~\ref{thm:woodroofe}
are all strong obstructions, existing in all dimensions $d\ge 1$.

For dimensions $\le 2$, by Theorem~\ref{thm:ots},
the strong obstructions to shellability are
the $1$-dimensional obstruction to shellability $2K_2$
and the simplicial complexes obtained by
adding 
zero or more
edges to (1a), (1b), (1c), (4a), (4b) and (4c)
of Figure~\ref{fig:obstructions}.
One advantage of considering strong obstructions instead of original
obstructions is that the number of strong obstructions is smaller than
that of original obstructions.
Another advantage is that the proof of identifying strong obstructions
is easier.

In the study of obstructions, 
some questions can be treated by discussions on strong obstructions.
In this section we discuss the following two questions on obstructions.

\begin{question} \label{question:finiteness}
For any fixed dimension, are there only finitely many obstructions to
shellability/partitionability/sequential Cohen-Macaulayness?
\end{question}

\begin{question} \label{question:difference}
Do the set of obstructions to shellability, 
that to partitionability
and that to sequential Cohen-Macaulayness
coincide for all dimensions?
\end{question}

Question~\ref{question:finiteness} is already asked by Wachs in \cite{W:2000}
for shellability.
(It is affirmatively conjectured in \cite[Conjecture 3.1.15]{W:2007}.)
Question~\ref{question:difference} is a natural question from 
our results of Theorem~\ref{thm:otp-otCM} for dimensions $\le 2$
and Theorem~\ref{thm:flag} for flag complexes.

In the following, we show these questions can be treated via 
discussions of the same questions for strong obstructions.

\begin{thm} \label{thm:s-obstruction1}
Let $\cal P$ be a property of simplicial complexes that is link-preserving.
If the number of strong obstructions to $\cal P$ 
of dimensions $\le d$ is finite, then
the number of obstructions to $\cal P$ 
of dimensions $\le d$ is also finite.
\end{thm}
\begin{proof}
Let $\Gamma$ be an obstruction to $\cal P$ but not a strong obstruction
to $\cal P$.
We first observe there is a vertex $v$ of $\Gamma$ such that 
$\link_\Gamma(\{v\})$ does not satisfy $\cal P$ as follows.
Since $\Gamma$ is an obstruction but not a strong obstruction to $\cal P$,
there is a face $\tau\neq\emptyset$ such that 
$\link_\Gamma(\tau)$ does not satisfy $\cal P$.
When $\lvert \tau\rvert=1$,
taking $\tau$ as the vertex $\{v\}$ satisfies the condition.
In the case $\lvert \tau\rvert \ge 2$,
let $v$ be one of the vertices of $\tau$. 
Then $\link_\Gamma(\{v\})$ does not satisfy $\cal P$,
since $\link_\Gamma(\tau) = \link_{\link_{\Gamma}(\{v\})}(\tau\sm\{v\})$
does not satisfy $\cal P$ and $\cal P$ is link-preserving.

For such a vertex $v$,
$\link_\Gamma(\{v\})$ is in fact an obstruction to $\cal P$.
For each $W\subsetneq V(\link_\Gamma(\{v\}))$,
$\link_\Gamma(\{v\})[W] = \link_{\Gamma[W]}(\{v\})$.
Since $\Gamma[W]$ satisfies $\cal P$,
$\link_\Gamma(\{v\})[W]$ satisfies $\cal P$.

Here we have $V(\Gamma) = \{v\} \cup V(\link_\Gamma(\{v\}))$,
since, if there is a vertex $w\notin \{v\} \cup V(\link_\Gamma(\{v\}))$, 
we have $\link_{\Gamma\sm \{w\}}(\{v\}) = \link_\Gamma(\{v\})$
and this implies that $\Gamma\sm \{w\}$ does not satisfy $\cal P$,
contradicting the fact $\Gamma$ is an obstruction to $\cal P$.

Now we show the statement of the theorem by induction on $d$.
For $d=0$, 
the number of obstructions to $\cal P$ is one (the case
that $\emptyset$ satisfies $\cal P$ and a singleton does not satisfy $\cal P$) 
or zero (otherwise), thus trivially finite. 
This constitute the base case of the induction.
Assume that the number of obstructions to $\cal P$ is finite
for each dimension $\le d-1$.
To show that the number of obstructions to $\cal P$ of dimension $d$ is finite,
we only have to show that the number of obstructions to $\cal P$ of dimension
$d$ that are not strong obstructions is finite.
The assumption that the number of obstructions to $\cal P$ of dimension $\le d-1$
is finite implies that there is a constant $C_{d-1} < \infty$ such that
$\lvert V(\Delta)\rvert \le C_{d-1}$ 
holds for each obstruction $\Delta$ to $\cal P$ of dimension $\le d-1$.
Now, for each obstruction $\Gamma$ to $\cal P$ of dimension $d$ 
that is not a strong obstruction, we have 
$V(\Gamma) = \{v\} \cup V(\link_\Gamma(\{v\}))$ and
$\link_\Gamma(\{v\})$ is an obstruction to $\cal P$ of dimension $\le d-1$.
This implies that $\lvert V(\Gamma)\rvert \le C_{d-1}+1$.
Thus the number of obstructions to $\cal P$ of dimension $d$
that are not strong obstructions is finite.
\end{proof}

\begin{thm} \label{thm:s-obstruction2}
Let $\cal P$ and $\cal Q$ be properties of simplicial complexes 
that are link-preserving.
Assume that $\cal P$ implies $\cal Q$.
If the set of strong obstructions to $\cal P$ and 
that of strong obstructions to $\cal Q$ coincide for dimensions $\le d$,
then the set of obstructions to $\cal P$ and 
that of obstructions to $\cal Q$ coincide for dimensions $\le d$.
\end{thm}

\begin{proof}
Let $\Gamma$ be an obstruction to $\cal P$ of dimension $\le d$.
If it is a strong obstruction to $\cal P$, 
then it is a strong obstruction to $\cal Q$ by assumption.
Especially, it does not satisfy $\cal Q$.

If $\Gamma$ is not a strong obstruction to $\cal P$,
there is a face $\tau\in \Gamma$ with $\tau\neq\emptyset$ 
such that $\link_\Gamma(\tau)$ does not
satisfy $\cal P$.
Take $\tau$ to be a maximal face with this property.
Then $\link_\Gamma(\tau)$ is a strong obstruction to $\cal P$:
$\link_\Gamma(\tau)[W]$ satisfies $\cal P$
for each $W\subsetneq V(\link_\Gamma(\{v\}))$
because $\link_\Gamma(\tau)[W] = \link_{\Gamma[W]}(\tau)$, 
and $\link_{\link_\Gamma(\tau)}(\eta)$ satisfies $\cal P$
for each $\eta\in \link_\Gamma(\tau)$ with $\eta\neq\emptyset$
because $\link_{\link_\Gamma(\tau)}(\eta)=\link_\Gamma(\tau\cup\eta)$
with $\tau \subsetneq \tau\cup\eta \in \Gamma$.
Since $\dim \link_\Gamma(\tau) < \dim\Gamma = d$, 
$\link_\Gamma(\tau)$ is a strong obstruction to $\cal Q$ by assumption,
thus $\link_\Gamma(\tau)$ does not satisfy $\cal Q$.
This implies that $\Gamma$ does not satisfy $\cal Q$.

Thus every obstruction to $\cal P$ does not satisfy $\cal Q$.
By Proposition~\ref{prop:difference} (by letting $\cal X$ be the class of
simplicial complexes of dimension $\le d$),
the set of obstructions to $\cal P$ equals to the set of obstructions
to $\cal Q$ in dimensions $\le d$.
\end{proof}

In closing this section, we add a proposition 
similar to Proposition~\ref{prop:difference},
which will be a convenient tool 
for a future discussion of the difference between 
classes of strong obstructions.
The proof can be done just in the same way
as Proposition~\ref{prop:difference}.

\begin{prop} \label{prop:difference-s}
Let $\cal P$ and $\cal Q$ be properties of simplicial complexes
and assume that $\cal P$ implies $\cal Q$.
In a class $\cal X$ of simplicial complexes that is closed under
restriction and link operation,
the following are equivalent.
\begin{enumerate}
\setlength{\itemsep}{0pt}
\item
In the class $\cal X$, the set of strong obstructions to $\cal P$
is different from the set of strong obstructions to $\cal Q$.
\item
There exists a strong obstruction to $\cal P$ in $\cal X$
which is not a strong obstruction to $\cal Q$.
\item
There exists a strong obstruction to $\cal P$ in $\cal X$
which satisfies $\cal Q$.
\end{enumerate}
\end{prop}

\medskip
\begin{remark}
Woodroofe~\cite{W:2009b} studies strong obstructions under the name
``dc-obstructions'', and enumerated those with up to 6 vertices 
with a help of computer.
\end{remark}

\section*{Acknowledgment}
The authors express their sincere gratitude to the anonymous referee
whose many valuable comments improved this manuscript.
During this work, the first author was partly supported by
KAKENHI \#18310101, \#19510137, \#20560052 and \#21740062.

\end{document}